\documentclass[12pt,a4paper]{article}
\usepackage{amsmath}
\usepackage{amssymb}
\usepackage{t1enc}
\usepackage[latin1]{inputenc}
\usepackage[english]{babel}
\usepackage{amsfonts}
\usepackage{latexsym}
\usepackage{enumerate}
\usepackage{bm}
\newtheorem{theorem}{Theorem}[section]
\newtheorem{prop}[theorem]{Proposition}
\newtheorem{lemma}[theorem]{Lemma}

\newtheorem{conj}[theorem]{Conjecture}

\newtheorem{rem}[theorem]{Remark}
\newtheorem{prob}[theorem]{Problem}
\newtheorem{example}[theorem]{Example}

\title{The maximum product of sizes of cross-intersecting families}
\author{Peter Borg\\[5mm]
Department of Mathematics, University of Malta, Malta\\
\texttt{peter.borg@um.edu.mt}}
\date{}

\begin{document}
\maketitle

\begin{abstract}
%
We say that a set $A$ \emph{$t$-intersects} a set $B$ if $A$ and $B$ have at least $t$ common elements. Two families $\mathcal{A}$ and $\mathcal{B}$ of sets are said to be \emph{cross-$t$-intersecting} if each set in $\mathcal{A}$ $t$-intersects each set in $\mathcal{B}$. A subfamily $\mathcal{S}$ of a family $\mathcal{F}$ is called a \emph{$t$-star of $\mathcal{F}$} if the sets in $\mathcal{S}$ have $t$ common elements. Let $l(\mathcal{F},t)$ denote the size of a largest $t$-star of $\mathcal{F}$. We call $\mathcal{F}$ a \emph{$(\leq r)$-family} if each set in $\mathcal{F}$ has at most $r$ elements. We determine 
a function $c : \mathbb{N}^3 \rightarrow \mathbb{N}$ such that the following holds. If $\mathcal{A}$ is a subfamily of a $(\leq r)$-family $\mathcal{F}$ with $l(\mathcal{F},t) \geq c(r,s,t)l(\mathcal{F},t+1)$, $\mathcal{B}$ is a subfamily of a $(\leq s)$-family $\mathcal{G}$ with $l(\mathcal{G},t) \geq c(r,s,t)l(\mathcal{G},t+1)$, and $\mathcal{A}$ and $\mathcal{B}$ are cross-$t$-intersecting, then $|\mathcal{A}||\mathcal{B}| \leq l(\mathcal{F},t)l(\mathcal{G},t)$. Some known results follow from this, and we 
identify several natural classes of families for which the bound is attained. 
\end{abstract}

\section{Introduction}

Unless otherwise stated, we shall use small letters such as $x$ to
denote non-negative integers or set elements or functions,
capital letters such as $X$ to denote sets, and calligraphic
letters such as $\mathcal{F}$ to denote \emph{families}
(that is, sets whose members are sets themselves). The set $\{1, 2, \dots\}$ of all positive integers is denoted by $\mathbb{N}$. 
For any $m,n \in \mathbb{N}$ with $m < n$, the set $\{i \in \mathbb{N} \colon m \leq i \leq n\}$ is denoted by $[m,n]$. We abbreviate $[1,n]$ to $[n]$. It is to be assumed that arbitrary sets and families are \emph{finite}. We call a set $A$ an \emph{$r$-element set}, or simply an \emph{$r$-set}, if its size $|A|$ is $r$. For a set $X$, the \emph{power set of $X$} (that is, the family of all subsets of $X$) is denoted by $2^X$, and the family of all $r$-element subsets of $X$ is denoted by ${X \choose r}$. 


We say that a set $A$ \emph{$t$-intersects} a set $B$ if $A$ and $B$ contain at least $t$ common elements. A family $\mathcal{A}$ of sets is said to be \emph{$t$-intersecting} if every two sets in $\mathcal{A}$ $t$-intersect. A $1$-intersecting family is also simply called an \emph{intersecting family}. 

For a family $\mathcal{F}$ and a set $T$, we denote the family $\{F \in \mathcal{F} \colon T \subseteq F\}$ by $\mathcal{F}(T)$. 
We call $\mathcal{F}(T)$ a \emph{$t$-star of $\mathcal{F}$} if $|T| = t$. 
A $t$-star of a family is the simplest example of a $t$-intersecting subfamily. We denote the size of a largest $t$-star of $\mathcal{F}$ by $l(\mathcal{F},t)$. We denote the set of largest $t$-stars of $\mathcal{F}$ by ${\rm L}(\mathcal{F},t)$. We say that $\mathcal{F}$ has the \emph{$t$-star property} if at least one $t$-star of $\mathcal{F}$ is a largest $t$-intersecting subfamily of $\mathcal{F}$.

One of the most popular endeavours in extremal set theory is that of determining the size or the structure of a largest $t$-intersecting subfamily of a given family $\mathcal{F}$. This originated in \cite{EKR}, which features the classical Erd\H os-Ko-Rado (EKR) Theorem. The EKR Theorem says that, for $1 \leq t \leq r$, there exists an integer $n_0(r,t)$ such that, for every $n \geq n_0(r,t)$, the size of a largest $t$-intersecting subfamily of ${[n] \choose r}$ is ${n-t \choose r-t}$, meaning that ${[n] \choose r}$ has the $t$-star property. It was also shown in \cite{EKR} that the smallest possible value of $n_0(r,1)$ is $2r$, and two of the various proofs of this fact (see \cite{Kat,K,D}) are particularly short and beautiful: Katona's \cite{K}, introducing the elegant cycle method, and Daykin's \cite{D}, using the Kruskal-Katona Theorem \cite{Kr,Ka}. If $n/2 < r < n$, then ${[n] \choose r}$ itself is intersecting.
A sequence of results \cite{EKR,F_t1,W,FF,AK1} culminated in the complete solution of the problem for $t$-intersecting subfamilies of ${[n] \choose r}$. The solution confirmed a conjecture of Frankl
\cite{F_t1} and 
particularly tells us that ${[n] \choose r}$ has the $t$-star property if and only if $n \geq (t+1)(r-t+1)$ \cite{F_t1,W}. The same $t$-intersection problem for $2^{[n]}$ was solved by Katona \cite{Kat}. 
These are among the most prominent results in extremal set theory. 
The EKR Theorem inspired a wealth of results that establish how large a system of sets can be under certain intersection conditions; see \cite{DF,F,F2,HST,HT,Borg7}.


Two families $\mathcal{A}$ and $\mathcal{B}$ are said to be \emph{cross-$t$-intersecting} if each set in $\mathcal{A}$ $t$-intersects each set in $\mathcal{B}$. More generally, $k$ families $\mathcal{A}_1, \dots, \mathcal{A}_k$ are said to be \emph{cross-$t$-intersecting} if for every $i$ and $j$ in $[k]$ with $i \neq j$, each set in $\mathcal{A}_i$ $t$-intersects each set in $\mathcal{A}_j$. Cross-$1$-intersecting families are also simply called \emph{cross-intersecting families}.

For $t$-intersecting subfamilies of a given family $\mathcal{F}$,
the natural question to ask is how large they can be. For
cross-$t$-intersecting families, two natural parameters arise: the
sum and the product of sizes of the cross-$t$-intersecting
families (note that the product of sizes of $k$ families
$\mathcal{A}_1, \dots, \mathcal{A}_k$ is the number of $k$-tuples
$(A_1, \dots, A_k)$ such that $A_i \in \mathcal{A}_i$ for each $i
\in [k]$). It is therefore natural to consider the problem of
maximising the sum or the product of sizes of $k$ cross-$t$-intersecting subfamilies (not necessarily distinct or non-empty) of a given family $\mathcal{F}$. The paper \cite{Borg8} analyses this problem in general, particularly reducing it to the problem of maximising the size of a $t$-intersecting subfamily of $\mathcal{F}$ for $k$ sufficiently large. Solutions have been obtained for various families (see \cite{Borg8}). 

Wang and Zhang \cite{WZ} solved the maximum sum problem for an important class of families that particularly includes ${[n] \choose r}$, using a striking combination of the method in \cite{Borg4,Borg3,Borg2,BL2,Borg5} and an important lemma that is found in \cite{AC,CK} and is referred to as the \emph{no-homomorphism lemma}. The solution for ${[n] \choose r}$ with $t=1$ had been obtained by Hilton \cite{H} and is the first result of this kind. 
For $2^{[n]}$, the maximum sum problem was solved \cite[Theorems~3.10, 4.1]{Borg8} via the result in \cite{WZ}, and the maximum product problem was settled in \cite{MT2} for the case where $k = 2$ or $n+t$ is even (see \cite[Section~5.2]{Borg8}, which features a conjecture for the case where $k > 2$ and $n+t$ is odd).

In this paper, we address the maximum product problem for the more general setting where each $\mathcal{A}_i$ is a subfamily of a family $\mathcal{F}_i$. This has been considered for a few special families \cite{Pyber,MT,Hirschorn,Borg11,BorgBLMS}, and, as we explain below, in many cases it is enough to solve the problem for $k = 2$ (see Lemma~\ref{prodgenlemma}).

The maximum product problem for ${[n] \choose r}$ was first addressed by Pyber \cite{Pyber}, who proved that, for $r, s, n \in \mathbb{N}$ such that either $r = s \leq n/2$ or $r < s$ and $n \geq 2s + r -2$, if $\mathcal{A} \subseteq {[n] \choose r}$ and $\mathcal{B} \subseteq {[n] \choose s}$ such that $\mathcal{A}$ and $\mathcal{B}$ are cross-intersecting, then $|\mathcal{A}||\mathcal{B}| \leq {n-1 \choose r-1}{n-1 \choose s-1}$. Subsequently, Matsumoto and Tokushige \cite{MT} proved this for $r \leq s \leq n/2$ (see also \cite{Bey}). 
For cross-$t$-intersecting subfamilies, we have the following. 

\begin{theorem}[\cite{Borg11}] \label{nchooser} For $1 \leq t \leq r \leq s$, there exists an integer $n_0(r,s,t)$ such that, for every $n \geq n_0(r,s,t)$, if $\mathcal{A} \subseteq {[n] \choose r}$, $\mathcal{B} \subseteq {[n] \choose s}$, and $\mathcal{A}$ and $\mathcal{B}$ are cross-$t$-intersecting, then $|\mathcal{A}||\mathcal{B}| \leq {n-t \choose r-t}{n-t \choose s-t}$, and equality holds if and only if $\mathcal{A} = \{A \in {[n] \choose r} \colon T \subseteq A\}$ and $\mathcal{B} = \{B \in {[n] \choose s} \colon T \subseteq B\}$ for some $T \in {[n] \choose t}$.
\end{theorem}
Hirschorn made a Frankl-type conjecture \cite[Conjecture~4]{Hirschorn} for any $r$, $s$, $t$ and $n$. A value of $n_0(r,s,t)$ that is close to best possible is established in \cite{Borg12}. The special case $r = s$ is treated in \cite{Tok1,Tok2,FLST}, which establish values of $n_0(r,r,t)$ that are also nearly optimal. 

Let $c : \mathbb{N}^3 \rightarrow \mathbb{N}$ such that, for $r,s,t \in \mathbb{N}$, $c(r,s,t) = \max \left\{ r{s \choose t}, s{r \choose t} \right\} + 1$ if $t \leq \min\{r,s\}$, and $c(r,s,t) = 1$ otherwise. Clearly, $c(r,s,t) = r{s \choose t}+1$ for $t \leq r \leq s$. 

The following is our main result, proved in Section~\ref{proofmain}.

\begin{theorem}\label{main} If $r,s,t \in \mathbb{N}$, $\mathcal{F}$ is a $(\leq r)$-family with $l(\mathcal{F},t) \geq c(r,s,t)l(\mathcal{F},t+1)$, $\mathcal{G}$ is a $(\leq s)$-family with $l(\mathcal{G},t) \geq c(r,s,t)l(\mathcal{G},t+1)$, and $\mathcal{A}$ and $\mathcal{B}$ are cross-$t$-intersecting families such that $\mathcal{A} \subseteq \mathcal{F}$ and $\mathcal{B} \subseteq \mathcal{G}$, then
$$|\mathcal{A}||\mathcal{B}| \leq l(\mathcal{F},t)l(\mathcal{G},t),$$
and equality holds if and only if $\mathcal{A} = \mathcal{F}(T) \in {\rm L}(\mathcal{F},t)$ and $\mathcal{B} = \mathcal{G}(T) \in {\rm L}(\mathcal{G},t)$ for some $t$-set $T$.
\end{theorem}
As we show in Section~\ref{specialfamilies}, this solves the problem for many natural families with a sufficiently large parameter depending on $r$, $s$ and $t$. 
For example, Theorem~\ref{main} yields Theorem~\ref{nchooser} by taking $n$ large enough so that ${n-t \choose r-t} \geq c(r,s,t){n-t-1 \choose r-t-1}$; see Section~\ref{lps}.

For $r, s, t \in \mathbb{N}$, let $\chi(r,s,t)$ be the smallest non-negative real number $a$ such that $|\mathcal{A}||\mathcal{B}| \leq l(\mathcal{F},t)l(\mathcal{G},t)$ for every $\mathcal{A}$, $\mathcal{B}$, $\mathcal{F}$ and $\mathcal{G}$ such that $\mathcal{F}$ is a $(\leq r)$-family with $l(\mathcal{F},t) \geq a \cdot l(\mathcal{F},t+1)$, $\mathcal{G}$ is a $(\leq s)$-family with $l(\mathcal{G},t) \geq a \cdot l(\mathcal{G},t+1)$, $\mathcal{A} \subseteq \mathcal{F}$, $\mathcal{B} \subseteq \mathcal{G}$, and $\mathcal{A}$ and $\mathcal{B}$ are cross-$t$-intersecting.

\begin{prob} What is the value of $\chi(r,s,t)$?
\end{prob}
By Theorem~\ref{main}, $\chi(r,s,t) \leq c(r,s,t)$.

In Theorem~\ref{main}, the case $\mathcal{F} = \mathcal{G}$ is of particular importance. First of all, it implies that $\mathcal{F}$ has the $t$-star property if $l(\mathcal{F},t) \geq c(r,r,t)l(\mathcal{F},t+1)$. 

\begin{theorem} \label{t-intversion} If $1 \leq t \leq r$ and $\mathcal{A}$ is a $t$-intersecting subfamily of a $(\leq r)$-family $\mathcal{F}$ with $l(\mathcal{F},t) \geq c(r,r,t)l(\mathcal{F},t+1)$, then 
\[|\mathcal{A}| \leq l(\mathcal{F},t),\]
and equality holds if and only if $\mathcal{A} \in {\rm L}(\mathcal{F},t)$. 
\end{theorem}
\textbf{Proof.} Let $\mathcal{G} = \mathcal{F}$ and $\mathcal{B} = \mathcal{A}$. Since $\mathcal{A}$ is $t$-intersecting, $\mathcal{A}$ and $\mathcal{B}$ are cross-$t$-intersecting. By Theorem~\ref{main}, the result follows.~\hfill{$\Box$} \\
 
Also note that in Theorem~\ref{main} with $\mathcal{F} = \mathcal{G}$, the bound is attained by taking $\mathcal{A} = \mathcal{B} \in {\rm L}(\mathcal{F},t)$; a generalization of this fact is given by Proposition~\ref{prop1}. As we show in Example~\ref{example1}, 
for $\mathcal{F} \neq \mathcal{G}$, it may be that the bound is not attained, and we may also have $\mathcal{A}$ and $\mathcal{B}$ for which no $t$-set $T$ satisfies $|\mathcal{A}||\mathcal{B}| \leq |\mathcal{F}(T)||\mathcal{G}(T)|$, no matter how large $\frac{l(\mathcal{F},t)}{l(\mathcal{F},t+1)}$ and $\frac{l(\mathcal{G},t)}{l(\mathcal{G},t+1)}$ are required to be. In view of this, we will now introduce further definitions. We will also generalise Theorem~\ref{main} to a result for $k$ cross-$t$-intersecting families.

\section{The cross-$t$-star property}

If $\mathcal{A}_1, \dots, \mathcal{A}_k$ are cross-$t$-intersecting families, then we say that the tuple $(\mathcal{A}_1, \dots, \mathcal{A}_k)$ is \emph{cross-$t$-intersecting}.

Let $\mathcal{F}_1, \dots, \mathcal{F}_k$ be families. We say that $(\mathcal{A}_1, \dots, \mathcal{A}_k)$ is \emph{below} $(\mathcal{F}_1, \dots, \mathcal{F}_k)$ if $\mathcal{A}_i \subseteq \mathcal{F}_i$ for each $i \in [k]$. We say that $(\mathcal{F}_1, \dots, \mathcal{F}_k)$ has the 
\begin {enumerate}[(a)]
\item \emph{cross-$t$-star property} if $\prod_{i=1}^k |\mathcal{A}_i| \leq \prod_{i=1}^k l(\mathcal{F}_i,t)$ for each cross-$t$-intersecting tuple $(\mathcal{A}_1, \dots, \mathcal{A}_k)$ below $(\mathcal{F}_1, \dots, \mathcal{F}_k)$.
\item \emph{strict cross-$t$-star property} if, for each cross-$t$-intersecting tuple $(\mathcal{A}_1, \dots, \mathcal{A}_k)$ below $(\mathcal{F}_1, \dots, \mathcal{F}_k)$, $\prod_{i=1}^k |\mathcal{A}_i| \leq \prod_{i=1}^k l(\mathcal{F}_i,t)$, and the inequality is strict if there exists no $t$-set $T$ such that, for each $i \in [k]$, $\mathcal{A}_i = \mathcal{F}_i(T)$.
\item \emph{strong cross-$t$-star property} if, for some $t$-set $T$, $\prod_{i=1}^k |\mathcal{A}_i| \leq \prod_{i=1}^k |\mathcal{F}_i(T)|$ for each cross-$t$-intersecting tuple $(\mathcal{A}_1, \dots, \mathcal{A}_k)$ below $(\mathcal{F}_1, \dots, \mathcal{F}_k)$.
\item \emph{extrastrong cross-$t$-star property} if there exists a $t$-set $T$ such that, for each cross-$t$-intersecting tuple $(\mathcal{A}_1, \dots, \mathcal{A}_k)$ below $(\mathcal{F}_1, \dots, \mathcal{F}_k)$, $\prod_{i=1}^k |\mathcal{A}_i| \leq \prod_{i=1}^k |\mathcal{F}_i(T)|$, and equality holds only if there exists a $t$-set $T'$ such that, for each $i \in [k]$, $\mathcal{A}_i = \mathcal{F}_i(T')$.
\end {enumerate}
Note that each of (b)--(d) implies (a), and (d) implies (a)--(c). As we demonstrate in Example~\ref{example1}, it may be that (b) holds, (c) does not hold, and hence (d) does not hold; clearly, this is the case only if $\prod_{i=1}^k |\mathcal{A}_i| < \prod_{i=1}^k l(\mathcal{F}_i,t)$ for each cross-$t$-intersecting tuple $(\mathcal{A}_1, \dots, \mathcal{A}_k)$ below $(\mathcal{F}_1, \dots, \mathcal{F}_k)$. 

\begin{example}\label{example1} \emph{Let $r_1, \dots, r_k, t \in \mathbb{N}$ with $k \geq 2$ and $t < r_1 \leq \dots \leq r_k$. Let $T_1, \dots, T_{k}, A_{1,1}, \dots, A_{1,q_1}, \dots, A_{k,1}, \dots, A_{k,q_k}$ be pairwise disjoint sets such that, for each $i \in [k]$, $|T_i| = t$ and $|A_{i,1}| = \dots = |A_{i,q_i}| = r_i-t$. Let $R_1, \dots, R_{k-1}$ be sets such that $|R_i| = r_i$ for each $i \in [k-1]$, $R_1 \subseteq \dots \subseteq R_{k-1}$, and $R_{k-1} \cap \bigcup_{i=1}^k\bigcup_{j=1}^{q_i} (T_i \cup A_{i,j}) = \emptyset$ (that is, no set $R_m$ intersects a set $T_i \cup A_{i,j}$). 
For each $i \in [k-1]$, let $\mathcal{F}_i = \{T_i \cup A_{i,1}, \dots, T_i \cup A_{i,q_i}, R_i\}$. Let $\mathcal{F}_k = \{T \cup A_{k,j} \colon T \in {R_1 \choose t}, j \in [q_k]\}$. For each $i \in [k]$, each set in $\mathcal{F}_i$ is of size $r_i$, and clearly $l(\mathcal{F}_i,t) = q_i$. For every $i,j \in [k]$ with $i < j$, a set $A$ in $\mathcal{F}_i$ $t$-intersects a set $B$ in $\mathcal{F}_j$ if and only if $A = R_i$ and either $j < k$ and $B = R_j$ or $j = k$ and $B \in \mathcal{F}_j$. 
Let $(\mathcal{A}_1, \dots, \mathcal{A}_k)$ be a cross-$t$-intersecting tuple below $(\mathcal{F}_1, \dots, \mathcal{F}_k)$ such that $\mathcal{A}_1, \dots, \mathcal{A}_k$ are non-empty (so that $\prod_{i=1}^k |\mathcal{A}_i| \neq 0$). Then $\mathcal{A}_i = \{R_i\}$ for each $i \in [k-1]$. Thus $\prod_{i=1}^k |\mathcal{A}_i| \leq |\mathcal{F}_k|$, and equality holds if and only if $(\mathcal{A}_1, \dots, \mathcal{A}_k) = (\{R_1\}, \dots, \{R_{k-1}\}, \mathcal{F}_k)$. Therefore, if $\prod_{i=1}^{k-1}q_i > {r_1 \choose t}$, then $\prod_{i=1}^k |\mathcal{A}_i| < \prod_{i=1}^k l(\mathcal{F}_i,t)$ (since $\prod_{i=1}^k |\mathcal{A}_i| \leq |\mathcal{F}_k| = {r_1 \choose t}q_k$ and $\prod_{i=1}^k l(\mathcal{F}_i,t) = \prod_{i=1}^k q_i$), meaning that $(\mathcal{F}_1, \dots, \mathcal{F}_k)$ has the strict cross-$t$-star property. Now let $T$ be a $t$-set such that $\prod_{i=1}^k |\mathcal{F}_i(T)| \neq 0$. Then $T \subseteq R_1$ and $\mathcal{F}_i(T) = \{R_i\}$ for each $i \in [k-1]$. Thus $\prod_{i=1}^k |\mathcal{F}_i(T)| = |\mathcal{F}_k(T)| < \prod_{i=1}^k |\mathcal{A}_i|$ if $(\mathcal{A}_1, \dots, \mathcal{A}_k) = (\{R_1\}, \dots, \{R_{k-1}\}, \mathcal{F}_k)$. Therefore, $(\mathcal{F}_1, \dots, \mathcal{F}_k)$ does not have the strong cross-$t$-star property.}
\end{example}
%

\begin{rem} \emph{By Example~\ref{example1}, for $1 \leq t < r \leq s$, there is no real number $a$ such that $(\mathcal{F},\mathcal{G})$ has the strong cross-$t$-star property for every $(\leq r)$-family $\mathcal{F}$ with $l(\mathcal{F},t) \geq a \cdot l(\mathcal{F},t+1)$ and every $(\leq s)$-family $\mathcal{G}$ with $l(\mathcal{G},t) \geq a \cdot l(\mathcal{G},t+1)$.}
\end{rem}



Of particular importance is the case $\mathcal{F}_1 = \dots = \mathcal{F}_k$.

\begin{prop}\label{prop1} (i) If $\mathcal{F}_1 = \dots = \mathcal{F}_k$ and $(\mathcal{F}_1, \dots, \mathcal{F}_k)$ has the cross-$t$-star property, then $(\mathcal{F}_1, \dots, \mathcal{F}_k)$ has the strong cross-$t$-star property. \\
(ii) If $\mathcal{F}_1 = \dots = \mathcal{F}_k$ and $(\mathcal{F}_1, \dots, \mathcal{F}_k)$ has the strict cross-$t$-star property, then $(\mathcal{F}_1, \dots, \mathcal{F}_k)$ has the extrastrong cross-$t$-star property.
\end{prop} 
\textbf{Proof.} Suppose $\mathcal{F}_1 = \dots = \mathcal{F}_k$. Let $(\mathcal{A}_1, \dots, \mathcal{A}_k)$ be a cross-$t$-intersecting tuple below $(\mathcal{F}_1, \dots, \mathcal{F}_k)$. Let $T$ be a $t$-set such that $|\mathcal{F}_1(T)| = l(\mathcal{F}_1,t)$. Since $\mathcal{F}_1 = \dots = \mathcal{F}_k$, $|\mathcal{F}_i(T)| = l(\mathcal{F}_i,t)$ for each $i \in [k]$. If $(\mathcal{F}_1, \dots, \mathcal{F}_k)$ has the cross-$t$-star property, then $\prod_{i=1}^k |\mathcal{A}_i| \leq \prod_{i=1}^k |\mathcal{F}_i(T)|$. If $(\mathcal{F}_1, \dots, \mathcal{F}_k)$ has the strict cross-$t$-star property and $\prod_{i=1}^k |\mathcal{A}_i| = \prod_{i=1}^k |\mathcal{F}_i(T)|$, then there exists a $t$-set $T'$ such that, for each $i \in [k]$, $\mathcal{A}_i = \mathcal{F}_i(T')$.~\hfill{$\Box$} 

\begin{prop}\label{prop2} The tuple $(\mathcal{F}_1, \dots, \mathcal{F}_k)$ has the extrastrong cross-$t$-star property if it has the strict cross-$t$-star property and there exists a $t$-set $T$ such that, for each $i \in [k]$, $\mathcal{F}_i(T) \in {\rm L}(\mathcal{F}_i,t)$.
\end{prop}
\textbf{Proof.} Let $(\mathcal{A}_1, \dots, \mathcal{A}_k)$ be a cross-$t$-intersecting tuple below $(\mathcal{F}_1, \dots, \mathcal{F}_k)$. Under the given conditions, $\prod_{i=1}^k |\mathcal{A}_i| \leq \prod_{i=1}^k l(\mathcal{F}_i,t) = \prod_{i=1}^k |\mathcal{F}_i(T)|$, and $\prod_{i=1}^k |\mathcal{A}_i| = \prod_{i=1}^k l(\mathcal{F}_i,t)$ only if there exists a $t$-set $T'$ such that, for each $i \in [k]$, $\mathcal{A}_i = \mathcal{F}_i(T')$.~\hfill{$\Box$}\\

By Theorem~\ref{main}, $(\mathcal{F},\mathcal{G})$ has the cross-$t$-star property if $l(\mathcal{F},t) \geq c(r,s,t)l(\mathcal{F},t+1)$ and $l(\mathcal{G},t) \geq c(r,s,t)l(\mathcal{G},t+1)$.

The following generalisation of \cite[Lemma~5.1]{Borg8} follows immediately from \cite[Lemma~5.2]{Borg8} and particularly tells us that the cross-$t$-star property is guaranteed for $k$ families if it is guaranteed for every two of them.

\begin{lemma} \label{prodgenlemma} If $2 \leq p \leq k$ and $\mathcal{F}_1, \dots, \mathcal{F}_k$ are families such that $(\mathcal{F}_{i_1}, \dots, \mathcal{F}_{i_p})$ has the cross-$t$-star property for each $p$-element subset $\{i_1, \dots, i_p\}$ of $[k]$, then $(\mathcal{F}_1, \dots, \mathcal{F}_k)$ has the cross-$t$-star property.
\end{lemma}


For example, 
Theorem~\ref{main} yields the following generalization.

\begin{theorem}\label{maingen} If $1 \leq t \leq r_1 \leq \dots \leq r_k$ and, for each $i \in [k]$, $\mathcal{F}_i$ is a $(\leq r_i)$-family with $l(\mathcal{F}_i,t) \geq c(r_{k-1},r_k,t)l(\mathcal{F}_i,t+1)$, then $(\mathcal{F}_1, \dots, \mathcal{F}_k)$ has the strict cross-$t$-star property.
\end{theorem}

We now start working towards the proofs of Theorems~\ref{main} and \ref{maingen}. Then, in Section~\ref{specialfamilies}, we apply the results above to several important families.

\section{Proof of the main result} \label{proofmain}

If a set $T$ $t$-intersects each set in a family $\mathcal{A}$, then we call $T$ a \emph{$t$-transversal of $\mathcal{A}$}. 

\begin{lemma}\label{main lemma} If $T$ is a $t$-transversal of a subfamily $\mathcal{A}$ of a family $\mathcal{F}$, then
\[|\mathcal{A}| \leq {|T| \choose t} l(\mathcal{F},t).\]
\end{lemma}
\textbf{Proof.} Let $\mathcal{T} = {T \choose t}$. Since $|A \cap T| \geq t$ for all $A \in \mathcal{A}$, we have
\begin{align} |\mathcal{A}| &= \left| \bigcup_{I \in \mathcal{T}}\mathcal{A}(I) \right| \leq \sum_{I \in
\mathcal{T}}|\mathcal{A}(I)| \leq \sum_{I \in \mathcal{T}}|\mathcal{F}(I)| \leq \sum_{I \in \mathcal{T}} l(\mathcal{F},t) = |\mathcal{T}| l(\mathcal{F},t), \nonumber
\end{align}
and hence the result.~\hfill{$\Box$}

\begin{lemma} \label{mainlemma2} If $T$ is a $t$-transversal of a subfamily $\mathcal{A}$ of a family $\mathcal{F}$, $X$ is a set of size $t$, $\mathcal{A} \subseteq \mathcal{F}(X)$, and $X \nsubseteq T$, then 
\[|\mathcal{A}| \leq |T \backslash X| l(\mathcal{F},t+1).\]
\end{lemma}
\textbf{Proof.} The result is trivial if $\mathcal{A} = \emptyset$. Suppose $\mathcal{A} \neq \emptyset$. For each $A \in \mathcal{A}$, we have 
\[t \leq |A \cap T| = |A \cap (T \cap X)| + |A \cap (T \backslash X)| = |T \cap X| + |A \cap (T \backslash X)| \leq t-1 + |A \cap (T \backslash X)|,\] 
and hence $|A \cap (T \backslash X)| \geq 1$. Together with $\mathcal{A} \subseteq \mathcal{F}(X)$, this gives us 
\begin{align} \mathcal{A} &\subseteq \{F \in \mathcal{F} \colon X \subseteq F, |F \cap (T \backslash X)| \geq 1\} \nonumber \\
&= \{F \in \mathcal{F} \colon X \cup \{y\}\subseteq F \mbox{ for some } y \in T \backslash X\} = \bigcup_{y \in T \backslash X} \mathcal{F}(X \cup \{y\}). \nonumber
\end{align}
Thus $|\mathcal{A}| \leq \sum_{y \in T \backslash X} |\mathcal{F}(X \cup \{y\})| \leq \sum_{y \in T \backslash X} l(\mathcal{F},t+1) = |T \backslash X| l(\mathcal{F},t+1)$.~\hfill{$\Box$}\\
\\
%
%

We can now prove Theorem~\ref{main}. We will call a $t$-intersecting family $\mathcal{A}$ \emph{trivial} if the sets in $\mathcal{A}$ have at least $t$ common elements.\\
\\
\textbf{Proof of Theorem~\ref{main}.} Suppose $|F| < t$ for each $F \in \mathcal{F}$. Then ${\rm L}(\mathcal{F},t) = \{\emptyset\}$, and hence $l(\mathcal{F},t) = 0 = l(\mathcal{F},t+1)$. Also, $|F \cap G| < t$ for each $F \in \mathcal{F}$ and each $G \in \mathcal{G}$. Thus, since $\mathcal{A}$ and $\mathcal{B}$ are cross-$t$-intersecting, one of $\mathcal{A}$ and $\mathcal{B}$ is empty, and hence $|\mathcal{A}||\mathcal{B}| = 0 = l(\mathcal{F},t)l(\mathcal{G},t)$. Similarly, $|\mathcal{A}||\mathcal{B}| = 0 = l(\mathcal{F},t)l(\mathcal{G},t)$ if $|G| < t$ for each $G \in \mathcal{G}$.

Now suppose that each of $\mathcal{F}$ and $\mathcal{G}$ has a set of size at least $t$. Then $r \geq t$, $s \geq t$, $l(\mathcal{F},t) \geq 1$ and $l(\mathcal{G},t) \geq 1$. 

If one of $\mathcal{A}$ and $\mathcal{B}$ is empty, then $|\mathcal{A}||\mathcal{B}| = 0 < l(\mathcal{F},t)l(\mathcal{G},t)$. 

Suppose $\mathcal{A} \neq \emptyset$ and $\mathcal{B} \neq \emptyset$. Since $\mathcal{A}$ and $\mathcal{B}$ are cross-$t$-intersecting, each set in $\mathcal{A}$ is a $t$-transversal of $\mathcal{B}$, and each set in $\mathcal{B}$ is a $t$-transversal of $\mathcal{A}$.\medskip

\textit{Case 1: $\mathcal{A}$ is not a trivial $t$-intersecting family, and $\mathcal{B}$ is not a trivial $t$-intersecting family.} Let $D \in \mathcal{B}$. For each $X \in {D \choose t}$, let $\mathcal{A}_X = \mathcal{A}(X)$. Since $|A \cap D| \geq t$ for each $A \in \mathcal{A}$, $\mathcal{A} = \bigcup_{X \in {D \choose t}} \mathcal{A}_X$. 

Consider any $X \in {D \choose t}$. Since $\mathcal{B}$ is not a trivial $t$-intersecting family, 
there exists $B \in \mathcal{B}$ such that $X \nsubseteq B$. Since $B$ is a $t$-transversal of $\mathcal{A}_X$, $|\mathcal{A}_X| \leq |B \backslash X| l(\mathcal{F},t+1) \leq s \cdot l(\mathcal{F},t+1)$ by Lemma~\ref{mainlemma2}. 

Therefore, we have
\begin{equation} |\mathcal{A}| = \left| \bigcup_{X \in {D \choose t}} \mathcal{A}_X \right| \leq \sum_{X \in {D \choose t}} |\mathcal{A}_X| \leq \sum_{X \in {D \choose t}} s l(\mathcal{F},t+1) = s {|D| \choose t} l(\mathcal{F},t+1), \nonumber
\end{equation}
and hence $|\mathcal{A}| \leq s {s \choose t} l(\mathcal{F},t+1)$. By a similar argument, $|\mathcal{B}| \leq r {r \choose t} l(\mathcal{G},t+1)$. It follows that
\begin{equation} |\mathcal{A}||\mathcal{B}| \leq s {s \choose t} l(\mathcal{F},t+1) r {r \choose t} l(\mathcal{G},t+1) \leq r {r \choose t} s {s \choose t} \frac{l(\mathcal{F},t)}{c(r,s,t)} \frac{l(\mathcal{G},t)}{c(r,s,t)} < l(\mathcal{F},t) l(\mathcal{G},t).  \nonumber
\end{equation}
%

\textit{Case 2: $\mathcal{A}$ is a trivial $t$-intersecting family, and $\mathcal{B}$ is not a trivial $t$-intersecting family.} We have $\mathcal{A} \subseteq \mathcal{F}(X)$ for some set $X$ of size $t$. Since $\mathcal{B}$ is not a trivial $t$-intersecting family, there exists $B \in \mathcal{B}$ such that $X \nsubseteq B$. By Lemma~\ref{mainlemma2}, $|\mathcal{A}| \leq |B \backslash X| l(\mathcal{F},t+1) \leq s \cdot l(\mathcal{F},t+1)$. Now let $C \in \mathcal{A}$. By Lemma~\ref{main lemma}, $|\mathcal{B}| \leq {|C| \choose t} l(\mathcal{G},t) \leq {r \choose t} l(\mathcal{G},t)$. Therefore, 
\begin{equation} |\mathcal{A}||\mathcal{B}| \leq s \cdot l(\mathcal{F},t+1){r \choose t} l(\mathcal{G},t) \leq s {r \choose t} \frac{l(\mathcal{F},t)}{c(r,s,t)} l(\mathcal{G},t) < l(\mathcal{F},t) l(\mathcal{G},t).  \nonumber
\end{equation}

\textit{Case 3: $\mathcal{A}$ is not a trivial $t$-intersecting family, and $\mathcal{B}$ is a trivial $t$-intersecting family.} The result follows by an argument similar to that for Case~2.\medskip

\textit{Case 4: $\mathcal{A}$ and $\mathcal{B}$ are trivial $t$-intersecting families.} 
Then $\mathcal{A} \subseteq \mathcal{F}(X)$ for some $t$-set $X$, and $\mathcal{B} \subseteq \mathcal{G}(Y)$ for some $t$-set $Y$. Thus $|\mathcal{A}| \leq l(\mathcal{F},t)$, $|\mathcal{B}| \leq l(\mathcal{G},t)$, and hence $|\mathcal{A}||\mathcal{B}| \leq l(\mathcal{F},t)l(\mathcal{G},t)$. Suppose $|\mathcal{A}||\mathcal{B}| = l(\mathcal{F},t)l(\mathcal{G},t)$. Then $|\mathcal{A}| = l(\mathcal{F},t)$ and $|\mathcal{B}| = l(\mathcal{G},t)$. Therefore, $\mathcal{A} = \mathcal{F}(X) \in {\rm L}(\mathcal{F},t)$ and $\mathcal{B} = \mathcal{G}(Y) \in {\rm L}(\mathcal{G},t)$. 

Suppose $X \neq Y$. Then $X \backslash Y \neq \emptyset$ since $|X| = |Y| = t$. Let $x \in X \backslash Y$. Suppose $x \in B$ for all $B \in \mathcal{B}$. Then $\mathcal{B} \subseteq \mathcal{G}(Y \cup \{x\})$, and hence $|\mathcal{B}| \leq l(\mathcal{G},t+1) \leq \frac{l(\mathcal{G},t)}{c(r,s,t)} < l(\mathcal{G},t)$, a contradiction. Thus $x \notin D$ for some $D \in \mathcal{B}$. Thus $X \nsubseteq D$. By Lemma~\ref{mainlemma2}, $|\mathcal{A}| \leq |D \backslash X| l(\mathcal{F},t+1) \leq s \frac{l(\mathcal{F},t)}{c(r,s,t)} < l(\mathcal{F},t)$, a contradiction.

Therefore, $X = Y$.~\hfill{$\Box$}\\
\\
\textbf{Proof of Theorem~\ref{maingen}.} For $k = 2$, the result is given by Theorem~\ref{main}. Consider $k \geq 3$. Let $(\mathcal{A}_1, \dots, \mathcal{A}_k)$ be a cross-$t$-intersecting tuple below $(\mathcal{F}_1, \dots, \mathcal{F}_k)$. Then, for every $i, j \in [k]$ with $i \neq j$, $\mathcal{A}_i$ and $\mathcal{A}_j$ are cross-$t$-intersecting, and, since $r_1 \leq \dots \leq r_k$, we have $c(r_i,r_j,t) \leq c(r_{k-1},r_k,t)$. By Theorem~\ref{main} and Lemma~\ref{prodgenlemma}, $\prod_{i=1}^k |\mathcal{A}_i| \leq \prod_{i=1}^k l(\mathcal{F}_i,t)$. Suppose equality holds.

Suppose $|\mathcal{A}_h| < l(\mathcal{F}_h,t)$ for some $h \in [k]$. By Theorem~\ref{main} and Lemma~\ref{prodgenlemma}, $\prod_{i \in [k] \backslash \{h\}} |\mathcal{A}_i| \leq \prod_{i \in [k] \backslash \{h\}} l(\mathcal{F}_i,t)$. Thus $\prod_{i=1}^k |\mathcal{A}_i| < \prod_{i=1}^k l(\mathcal{F}_i,t)$, a contradiction.

Therefore, $|\mathcal{A}_i| \geq l(\mathcal{F}_i,t)$ for each $i \in [k]$. Since $\prod_{i=1}^k |\mathcal{A}_i| = \prod_{i=1}^k l(\mathcal{F}_i,t)$, $|\mathcal{A}_i| = l(\mathcal{F}_i,t)$ for each $i \in [k]$. For each $i \in [2,k]$, we have $|\mathcal{A}_1||\mathcal{A}_i| = l(\mathcal{F}_1,t)l(\mathcal{F}_i,t)$, and hence, by Theorem~\ref{main}, there exists a $t$-set $T_{1,i}$ such that $\mathcal{A}_1 = \mathcal{F}_1(T_{1,i}) \in {\rm L}(\mathcal{F}_1,t)$ and $\mathcal{A}_i = \mathcal{F}_i(T_{1,i}) \in {\rm L}(\mathcal{F}_i,t)$. By the argument in Case~4 of the proof of Theorem~\ref{main}, $T_{1,i} = T_{1,2}$ for each $i \in [2,k]$.  
Thus $\mathcal{A}_i = \mathcal{F}_i(T_{1,2})$ for each $i \in [k]$.~\hfill{$\Box$}

\section{Classes of families} \label{specialfamilies}

In this section, we apply Theorem~\ref{maingen} to important classes of families. Thus, for each family $\mathcal{F}$, we need to obtain an upper bound for $\frac{l(\mathcal{F},t)}{l(\mathcal{F},t+1)}$ and compare it with $c(r,s,t)$. 

Much of the work done on the $t$-intersection problem for the families treated here is outlined in \cite{Borg7}. Much less is known about the product cross-$t$-intersection problem because it takes the $t$-intersection problem to a deeper level; most of the main results are outlined in \cite{Borg8}. 
We will show that Theorem~\ref{maingen} provides a solution for many of the most natural and mostly studied classes of families. For each class, Theorem~\ref{t-intversion} provides a solution for the $t$-intersection problem.

\subsection{Levels of power sets} \label{lps}

For a family $\mathcal{F}$ and a non-negative integer $r$, the family of all $r$-element sets in $\mathcal{F}$ is called the \emph{$r$-th level of $\mathcal{F}$}. For a set $X$, ${X \choose r}$ is the $r$-th level of $2^X$. 

Consider $\mathcal{F} = {[n] \choose p}$ with $1 \leq p \leq n$. Suppose $l(\mathcal{F},t+1) > 0$ for some $t \geq 1$. Then $p \geq t+1$. We have
\begin{equation} \frac{l(\mathcal{F},t)}{l(\mathcal{F},t+1)} = 
\frac{ {n-t \choose p-t} }{ {n-t-1 \choose p-t-1} } = \frac{n-t}{p-t}. \label{lps1}
\end{equation}
Therefore, $l(\mathcal{F},t) \geq c(r,s,t)l(\mathcal{F},t+1)$ if $n \geq (p-t)c(r,s,t) + t$.

%
%

The following is a generalization of Theorem~\ref{nchooser}.

\begin{theorem} \label{nchoosergen} If $1 \leq t \leq r_1 \leq \dots \leq r_k$ and, for each $i \in [k]$, $\mathcal{F}_i = {[n_i] \choose r_i}$ with $n_i \geq (r_i-t)c(r_{k-1},r_k,t) + t$, then $(\mathcal{F}_1, \dots, \mathcal{F}_k)$ has the extrastrong cross-$t$-star property.
\end{theorem}
\textbf{Proof.} By (\ref{lps1}), for each $i \in [k]$, $l(\mathcal{F}_i,t) \geq c(r_{k-1},r_k,t)l(\mathcal{F}_i,t+1)$ as $n_i \geq (r_i-t)c(r_{k-1},r_k,t) + t$. By Theorem~\ref{maingen}, $(\mathcal{F}_1, \dots, \mathcal{F}_k)$ has the strict cross-$t$-star property. Since $\mathcal{F}_i([t]) \in {\rm L}(\mathcal{F}_i,t)$ for each $i \in [k]$, the result follows by Proposition~\ref{prop2}.~\hfill{$\Box$}

\subsection{Families of integer sequences} 

For an $r$-element set $X = \{x_1, \dots, x_r\}$ 
and an integer $m \geq 1$, we define 
\[\mathcal{S}_{X,m} = \{\{(x_1,y_1), \dots, (x_r,y_r)\}
\colon y_1, \dots, y_r \in [m]\}.\]
Note that $\mathcal{S}_{X,m}$ is isomorphic to the set $[m]^r$, that is, the set of all sequences $(y_1, \dots, y_r)$ such that $y_i \in [m]$ for each $i \in [r]$. We take $\mathcal{S}_{\emptyset,m}$ to be $\emptyset$. With a slight abuse of notation, for a family $\mathcal{F}$, we define
\[\mathcal{S}_{\mathcal{F},m} = \bigcup_{F \in
\mathcal{F}}\mathcal{S}_{F,m}.\]
%

The $t$-intersection problem for $\mathcal{S}_{[n],m}$ was solved by Ahlswede and Khachatrian \cite{AK2} and by Frankl and Tokushige \cite{FT2}, and that for $\mathcal{S}_{\mathcal{F},m}$ is solved in \cite{Borg6} for $m$ sufficiently large. The product cross-$t$-intersection problem for $\mathcal{S}_{[n],m}$ was solved by Moon \cite{Moon} for $m \geq t+2$, and by Frankl et al.~\cite{FLST} and Pach and Tardos \cite{PT} for $m \geq t+1$. We solve the problem for $\mathcal{S}_{\mathcal{F},m}$ with $m$ sufficiently large depending only on $t$ and the size of a largest set in $\mathcal{F}$.

\begin{theorem} \label{intseqresult} If $1 \leq t \leq r_1 \leq \dots \leq r_k$ and, for each $i \in [k]$, $\mathcal{F}_i$ is a $(\leq r_i)$-family and $m_i \geq c(r_{k-1},r_k,t)$, then $(\mathcal{S}_{\mathcal{F}_1,m_1}, \dots, \mathcal{S}_{\mathcal{F}_k,m_k})$ has the strict cross-$t$-star property.
\end{theorem}
\begin{lemma} \label{intseqlemma} If $1 \leq t \leq r$ and $\mathcal{F}$ is a $(\leq r)$-family, then $l(\mathcal{S}_{\mathcal{F},m},t) \geq m \cdot l(\mathcal{S}_{\mathcal{F},m},t+1)$.
\end{lemma}
\textbf{Proof.} Suppose $l(\mathcal{S}_{\mathcal{F},m},t+1) > 0$. Then $r \geq t+1$. Let $\mathcal{A}$ be a $(t+1)$-star of $\mathcal{S}_{\mathcal{F},m}$ of size $l(\mathcal{S}_{\mathcal{F},m},t+1)$. Then $\mathcal{A} = \mathcal{S}_{\mathcal{F},m}(Z)$ for some $(t+1)$-element set $Z$. Let $\mathcal{G} = \{F \in \mathcal{F} \colon \mathcal{S}_{F,m}(Z) \neq \emptyset\}$. Let $T \in {Z \choose t}$. 
We have
\begin{align} l(\mathcal{S}_{\mathcal{F},m},t+1) &= |\mathcal{S}_{\mathcal{F},m}(Z)| = \sum_{F \in \mathcal{F}} |\mathcal{S}_{F,m}(Z)| = \sum_{F \in \mathcal{G}} |\mathcal{S}_{F,m}(Z)| = \sum_{F \in \mathcal{G}} m^{|F|-t-1} \nonumber \\
&= \frac{1}{m} \sum_{F \in \mathcal{G}} m^{|F|-t} = \frac{1}{m} \sum_{F \in \mathcal{G}} |\mathcal{S}_{F,m}(T)| \leq \frac{1}{m} \sum_{F \in \mathcal{F}} |\mathcal{S}_{F,m}(T)| \nonumber \\
&=\frac{1}{m} |\mathcal{S}_{\mathcal{F},m}(T)| \leq \frac{1}{m}l(\mathcal{S}_{\mathcal{F},m},t), \nonumber
\end{align}
and hence the result.~\hfill{$\Box$}\\
\\
\textbf{Proof of Theorem~\ref{intseqresult}.} For any $i \in [k]$, $l(\mathcal{S}_{\mathcal{F}_i,m_i},t) \geq c(r_{k-1},r_k,t) l(\mathcal{S}_{\mathcal{F}_i,m_i},t+1)$ by Lemma~\ref{intseqlemma} and the given condition $m_i \geq c(r_{k-1},r_k,t)$. The result follows by Theorem~\ref{maingen}.~\hfill{$\Box$}

\begin{theorem} \label{intseqresultcor} If $1 \leq t \leq r$, $\mathcal{F}$ is a $(\leq r)$-family, $m \geq c(r,r,t)$, and $\mathcal{F}_1 = \dots = \mathcal{F}_k = \mathcal{S}_{\mathcal{F},m}$, then $(\mathcal{F}_1, \dots, \mathcal{F}_k)$ has the extrastrong cross-$t$-star property.
\end{theorem}
\textbf{Proof.} The result follows by Theorem~\ref{intseqresult} and Proposition~\ref{prop1}.~\hfill{$\Box$} \\

We make the following conjecture, which is analogous to \cite[Conjecture~2.1]{Borg6}.

\begin{conj} \label{fisconj} For any $t \geq 1$, there
exists a positive integer $m_0(t)$ such that $(\mathcal{S}_{\mathcal{F},m}, \mathcal{S}_{\mathcal{F},m})$ has the strong cross-$t$-star property for any family $\mathcal{F}$ and any $m \geq m_0(t)$.
\end{conj}
We also conjecture that the smallest possible $m_0(t)$ is $t+1$, and that $(\mathcal{S}_{\mathcal{F},m}, \mathcal{S}_{\mathcal{F},m})$ has the extrastrong cross-$t$-star property if $m > t+1$. By Lemma~\ref{prodgenlemma} and Proposition~\ref{prop1}, this would imply a strengthening of Theorem~\ref{intseqresultcor}. The conjecture does not hold for $m < t+1$. Indeed, it can be checked that, if $m \leq t$, $n \geq t+2$, $\mathcal{F} = \{[n]\}$, and $\mathcal{A} = \mathcal{B} = \{A \in \mathcal{S}_{[n],m} \colon |A \cap \{(1,1), \dots, (t+2,1)\}| \geq t+1\}$, then $\mathcal{A}$ and $\mathcal{B}$ are cross-$t$-intersecting subfamilies of $\mathcal{S}_{\mathcal{F},m}$ and $|\mathcal{A}||\mathcal{B}| > (m^{n-t})^2 = (l(\mathcal{S}_{\mathcal{F},m}))^2$.

\subsection{Families of permutations} 

For an $r$-set $X = \{x_1, \dots, x_r\}$ and an integer $m \geq 1$, we define $\mathcal{S}_{X,m}^*$ to be the special subfamily of
$\mathcal{S}_{X,m}$ given by
$$\mathcal{S}_{X,m}^* = \left\{\{(x_1, y_1), \dots, (x_r, y_r)\}
\colon y_1, \dots, y_r \mbox{ are distinct elements of } [m] \right\}.$$
Note that $\mathcal{S}_{X,m}^* \neq \emptyset$ if and only if $r \leq m$. The family $\mathcal{S}_{X,m}^*$ can be interpreted as the set of permutations of sets in ${[m] \choose r}$; indeed, a member $\{(x_1, y_1), \dots, (x_r, y_r)\}$ of $\mathcal{S}_{X,m}^*$ corresponds uniquely to the permutation $(y_1, \dots, y_r)$ of the $r$-element subset $\{y_1, \dots, y_r\}$ of $[m]$. We take $\mathcal{S}_{\emptyset,m}$ to be $\emptyset$. With a slight abuse of notation, for a family $\mathcal{F}$, we define
$\mathcal{S}_{\mathcal{F},m}^*$ to be the special subfamily of
$\mathcal{S}_{\mathcal{F},m}$ given by
\[\mathcal{S}_{\mathcal{F},k}^* = \bigcup_{F \in
\mathcal{F}}\mathcal{S}_{F,k}^*.\]


In \cite{DF1}, Deza and Frankl established the $1$-star property of $\mathcal{S}_{[m],m}^*$ and conjectured that $\mathcal{S}_{[m],m}^*$ has the $t$-star property for $m$ sufficiently large depending on $t$. 
Ellis, Friedgut and Pilpel \cite{EFP} proved the conjecture together with the product cross-$t$-intersection version. The $t$-intersection problem for $\mathcal{S}_{\mathcal{F},m}^*$ is solved in \cite{Borg6} for $m$ sufficiently large depending only on $t$ and the size of a largest set in $\mathcal{F}$. For this case, we have the following analogous result for the product cross-$t$-intersection problem.

\begin{theorem} \label{permresult} If $1 \leq t \leq r_1 \leq \dots \leq r_k$ and, for each $i \in [k]$, $\mathcal{F}_i$ is a $(\leq r_i)$-family and $m_i \geq c(r_{k-1},r_k,t)+t$, then $(\mathcal{S}_{\mathcal{F}_1,m_1}^*, \dots, \mathcal{S}_{\mathcal{F}_k,m_k}^*)$ has the strict cross-$t$-star property.
\end{theorem}
Similarly to Theorem~\ref{intseqresult}, this follows from the fact that if 
$\mathcal{S}_{F,m}^*(Z) \neq \emptyset$ for some $(t+1)$-element $Z$, then $|\mathcal{S}_{F,m}^*(T)| = \frac{(m-t)!}{(m-|F|)!} = (m-t)|\mathcal{S}_{F,m}^*(Z)|$ for any $T \in {Z \choose t}$.

\begin{theorem} \label{permresultcor} If $1 \leq t \leq r$, $\mathcal{F}$ is a $(\leq r)$-family, $m \geq c(r,r,t) + t$, and $\mathcal{F}_1 = \dots = \mathcal{F}_k = \mathcal{S}_{\mathcal{F},m}^*$, then $(\mathcal{F}_1, \dots, \mathcal{F}_k)$ has the extrastrong cross-$t$-star property.
\end{theorem}
\textbf{Proof.} The result follows by Theorem~\ref{permresult} and Proposition~\ref{prop1}.~\hfill{$\Box$} \\

We make the following conjecture, which is analogous to \cite[Conjecture~2.4]{Borg6}.

\begin{conj} \label{fisconj} For any $t \geq 1$, there
exists a positive integer $m_0^*(t)$ such that $(\mathcal{S}_{\mathcal{F},m}^*, \mathcal{S}_{\mathcal{F},m}^*)$ has the extrastrong cross-$t$-star property for any family $\mathcal{F}$ and any $m \geq m_0^*(t)$.
\end{conj}
By Lemma~\ref{prodgenlemma} and Proposition~\ref{prop1}, this would imply a strengthening of Theorem~\ref{permresultcor}.

\subsection{Families of multisets}

A \emph{multiset} is a collection $A$ of objects such that each object possibly appears more than once in $A$. Thus the difference between a multiset and a set is that a multiset may have repetitions of its members. The \emph{multiplicity} of a member $a$ of a multiset $A$ is the number of instances of $a$ in $A$, and is denoted by $m_A(a)$. If $a_1, \dots, a_r$ are the distinct members of a multiset $A$, then we can represent $A$ uniquely by the set $\{(a_i,j) \colon i \in [r], j \in [m_A(a)]\}$, which we denote by $S_A$. Let $M_{n,r}$ denote the set of all multisets $A$ such that the members of $A$ are in $[n]$ and amount to $r$ with repetitions included. An elementary counting result is that \[|M_{n,r}| = {n+r-1 \choose r}.\] Let $\mathcal{M}_{n,r}$ denote the family $\{S_A \colon A \in M_{n,r}\}$. Note that two multisets $A$ and $B$ have exactly $q$ common members (with repetitions included) if and only if $|S_A \cap S_B| = q$.

The $t$-intersection problem for $\mathcal{M}_{n,r}$ was solved by Meagher and Purdy \cite{MP} for $t = 1$, and by F\"{u}redi, Gerbner and Vizer \cite{FGV} for $n \geq 2r-t$. Here we solve the product cross-$t$-intersection problem for $n$ sufficiently large depending on $r$ and $t$.

Consider $\mathcal{F} = \mathcal{M}_{n,p}$. Suppose $l(\mathcal{F},t+1) > 0$ for some $t \geq 1$. Then $p \geq t+1$. We have
\begin{equation} \frac{l(\mathcal{F},t)}{l(\mathcal{F},t+1)} = 
\frac{ {n+p-t-1 \choose p-t} }{ {n+p-t-2 \choose p-t-1} } = \frac{n+p-t-1}{p-t}. \label{fm1}
\end{equation}
Therefore, $l(\mathcal{F},t) \geq c(r,s,t)l(\mathcal{F},t+1)$ if $n \geq (p-t)c(r,s,t) - p + t + 1$.

\begin{theorem} \label{fmresult} If $1 \leq t \leq r_1 \leq \dots \leq r_k$ and, for each $i \in [k]$, $\mathcal{F}_i = \mathcal{M}_{n_i,r_i}$ with $n_i \geq (r_i-t)c(r_{k-1},r_k,t) - r_i + t + 1$, then $(\mathcal{F}_1, \dots, \mathcal{F}_k)$ has the extrastrong cross-$t$-star property.
\end{theorem}
\textbf{Proof.} For each $i \in [k]$, $l(\mathcal{F}_i,t) \geq c(r_{k-1},r_k,t)l(\mathcal{F}_i,t+1)$ by (\ref{fm1}) and the given condition $n_i \geq (r_i-t)c(r_{k-1},r_k,t) - r_i + t + 1$. By Theorem~\ref{maingen}, $(\mathcal{F}_1, \dots, \mathcal{F}_k)$ has the strict cross-$t$-star property. Let $T = \{(1,i) \colon i \in [t]\}$. Since $\mathcal{F}_i(T) \in {\rm L}(\mathcal{F}_i,t)$ for each $i \in [k]$, the result follows by Proposition~\ref{prop2}.~\hfill{$\Box$}

\subsection{Families of compositions}

If $a_1, a_2, \dots, a_r$ and $n$ are positive integers such that $n = a_1 + a_2 + \dots + a_r$, then the tuple $(a_1, a_2, \dots, a_r)$ is said to be a \emph{composition of $n$} of \emph{length $r$}. Let $C_{n,r}$ denote the set of all compositions of $n$ of length $r$. 
An elementary counting result is that \[|C_{n,r}| = {n-1 \choose n-r} = {n-1 \choose r-1}.\] We can represent a composition ${\bf a} = (a_1, \dots, a_r)$ uniquely by the set $\{(1,a_1), \dots, (r,a_r)\}$, which we denote by $S_{\bf a}$. Let $\mathcal{C}_{n,r}$ denote the family $\{S_{\bf a} \colon {\bf a} \in C_{n,r}\}$. 

We say that a composition ${\bf a} = (a_1, \dots, a_r)$ \emph{strongly $t$-intersects} a composition ${\bf b} = (b_1, \dots, b_s)$ if there exists a $t$-element subset $T$ of $[\min\{r,s\}]$ such that $a_i = b_i$ for each $i \in T$. Note that ${\bf a}$ \emph{strongly $t$-intersects} ${\bf b}$ if and only if $|S_{\bf a} \cap S_{\bf b}| \geq t$.

Ku and Wong \cite{KW1} solved the $t$-intersection problem for $\mathcal{C}_{n,r}$ with $n$ sufficiently large. In \cite{KW2}, they also proved Theorem~\ref{fcresult} below for sufficiently large values of $n_1, \dots, n_r$.

Consider $\mathcal{F} = \mathcal{C}_{n,p}$ with $t+1 < p \leq n$. It is straightforward that $\mathcal{F}(\{(i,1) \colon i \in [t]\})$ is a largest $t$-star of $\mathcal{C}_{n,p}$. We have
\begin{equation} \frac{l(\mathcal{F},t)}{l(\mathcal{F},t+1)} = 
\frac{ {n-t-1 \choose p-t-1} }{ {n-t-2 \choose p-t-2} } = \frac{n-t-1}{p-t-1}. \label{fc1}
\end{equation}
Therefore, $l(\mathcal{F},t) \geq c(r,s,t)l(\mathcal{F},t+1)$ if $n \geq (p-t-1)c(r,s,t) + t + 1$.

\begin{theorem} \label{fcresult} If $2 \leq t+1 < r_1 \leq \dots \leq r_k$ and, for each $i \in [k]$, $\mathcal{F}_i = \mathcal{C}_{n_i,r_i}$ with $n_i \geq (r_i-t-1)c(r_{k-1},r_k,t) + t + 1$, then $(\mathcal{F}_1, \dots, \mathcal{F}_k)$ has the extrastrong cross-$t$-star property.
\end{theorem}
\textbf{Proof.} For each $i \in [k]$, $l(\mathcal{F}_i,t) \geq c(r_{k-1},r_k,t)l(\mathcal{F}_i,t+1)$ by (\ref{fc1}) and the given condition $n_i \geq (r_i-t-1)c(r_{k-1},r_k,t) + t + 1$. By Theorem~\ref{maingen}, $(\mathcal{F}_1, \dots, \mathcal{F}_k)$ has the strict cross-$t$-star property. Let $T = \{(i,1) \colon i \in [t]\}$. Since $\mathcal{F}_i(T) \in {\rm L}(\mathcal{F}_i,t)$ for each $i \in [k]$, the result follows by Proposition~\ref{prop2}.~\hfill{$\Box$}

\subsection{Families of set partitions} 

If $X_1, X_2, \dots, X_r$ are pairwise disjoint non-empty sets and $X = \bigcup_{i=1}^r X_i$, then the set $\{X_1, X_2, \dots, X_r\}$ is called a \emph{partition of $X$} of \emph{length $r$}, and $X_1, X_2, \dots, X_r$ are called the \emph{parts} of the partition. Let $\mathsf{P}_{n,r}$ denote the family of all partitions of $[n]$ of length $r$, and let $s_{n,r} = |\mathsf{P}_{n,r}|$.  Trivially, $s_{n,1} = 1 = s_{n,n}$. An elementary result is that
\begin{equation} s_{n,r} = s_{n-1,r-1} + rs_{n-1,r} \quad \mbox{if } 2 \leq r \leq n-1. \nonumber
\end{equation}
It follows that 
\begin{equation} s_{m,r} \leq s_{n,r} \quad \mbox{if } 1 \leq m \leq n. \label{fsp1}
\end{equation}

\begin{lemma} \label{fsplemma1} If $1 < r < n$, then $s_{n,r} \geq \frac{n-1}{r-1}s_{n-1,r-1}$.
%
%
\end{lemma}
\textbf{Proof.} Consider any $X \in \mathsf{P}_{n-1,r-1}$. For any $i \in [n-1]$, let $A_i$ be the part of $X$ that contains $i$, and let $X_i$ be the member of $\mathsf{P}_{n,r}$ obtained by replacing $i$ by $n$ in $A_i$, and adding $\{i\}$ as a part; that is, $X_i = (X \backslash \{A_i\}) \cup \{(A_i \backslash \{i\}) \cup \{n\}\} \cup \{\{i\}\}$. For any $Y \in \mathsf{P}_{n,r}$, let $f(X_i,Y) = 1$ if $X_i = Y$, and let $f(X_i,Y) = 0$ if $X_i \neq Y$. If $Y$ has no parts of size $1$, then $f(X_i,Y) = 0$. Suppose that $\{y_1\}, \dots, \{y_p\}$ are the distinct parts of $Y$ of size $1$. Since $n > r$, $p \leq r-1$. Let $B$ be the part of $Y$ that contains $n$. Then $f(X_i,Y) = 1$ if and only if $i \in \{y_1, \dots, y_p\}$ and $X = (Y \backslash \{B, \{i\}\}) \cup \{(B \backslash \{n\}) \cup \{i\}\}$. 

Therefore, we have
\begin{align} (n-1) s_{n-1,r-1} &= \sum_{X \in \mathsf{P}_{n-1,r-1}} \sum_{i = 1}^{n-1} 1 = \sum_{X \in \mathsf{P}_{n-1,r-1}} \sum_{i = 1}^{n-1} \sum_{Y \in \mathsf{P}_{n,r}} f(X_i,Y) \nonumber \\
&= \sum_{Y \in \mathsf{P}_{n,r}} \sum_{X \in \mathsf{P}_{n-1,r-1}} \sum_{i = 1}^{n-1} f(X_i,Y) \leq \sum_{Y \in \mathsf{P}_{n,r}} (r-1) = (r-1) s_{n,r}, \nonumber
\end{align}
and hence the result.~\hfill{$\Box$}\\

Erd\H{o}s and Sz\'{e}kely \cite{ES} solved the $t$-intersection problem for $\mathsf{P}_{n,r}$ with $n$ sufficiently large (see \cite{KR} for a related result). Using the results above, we prove the following cross-$t$-intersection result.

\begin{theorem} \label{fspresult} If $2 \leq t+1 < r_1 \leq \dots \leq r_k$ and, for each $i \in [k]$, $\mathcal{F}_i = \mathsf{P}_{n_i,r_i}$ with $n_i \geq (r_i-t-1)c(r_{k-1},r_k,t) + t + 1$, then $(\mathcal{F}_1, \dots, \mathcal{F}_k)$ has the extrastrong cross-$t$-star property.
\end{theorem}

\begin{lemma} \label{fsplemma2} If $1 \leq t < r \leq n$, then $l(\mathsf{P}_{n,r},t) = s_{n-t,r-t}$.
\end{lemma}
\textbf{Proof.} Let $T = \{\{i\} \colon i \in [t]\}$. We have $l(\mathsf{P}_{n,r},t) \geq |\mathsf{P}_{n,r}(T)| = s_{n-t,r-t}$. Let $\mathcal{A}$ be a largest $t$-star of $\mathsf{P}_{n,r}$. There exist $t$ pairwise disjoint non-empty subsets $X_1, \dots, X_t$ of $[n]$ such that $\mathcal{A} = \mathsf{P}_{n,r}(\{X_1, \dots, X_t\})$. Thus $l(\mathsf{P}_{n,r},t) = |\mathcal{A}| = s_{n',r-t}$, where $n' = n - \sum_{i=1}^t |X_i| \leq n-t$. By (\ref{fsp1}), $l(\mathsf{P}_{n,r},t) \leq s_{n-t,r-t}$. Since $l(\mathsf{P}_{n,r},t) \geq s_{n-t,r-t}$, the result follows.~\hfill{$\Box$}

\begin{lemma} \label{fsplemma3} If $2 \leq t+1 < r < n$, then $l(\mathsf{P}_{n,r},t) \geq \frac{n-t-1}{r-t-1}l(\mathsf{P}_{n,r},t+1)$.
\end{lemma}
\textbf{Proof.} 
%
By Lemma~\ref{fsplemma2}, $l(\mathsf{P}_{n,r},t) = s_{n-t,r-t}$ and $l(\mathsf{P}_{n,r},t+1) = s_{n-t-1,r-t-1}$. Thus, by Lemma~\ref{fsplemma1}, $l(\mathsf{P}_{n,r},t) \geq \frac{n-t-1}{r-t-1} l(\mathsf{P}_{n,r},t+1)$.~\hfill{$\Box$}\\
\\ 
\textbf{Proof of Theorem~\ref{fspresult}.} For each $i \in [k]$, $l(\mathcal{F}_i,t) \geq c(r_{k-1},r_k,t)l(\mathcal{F}_i,t+1)$ by Lemma~\ref{fsplemma3} and the given condition $n_i \geq (r_i-t-1)c(r_{k-1},r_k,t) + t + 1$. By Theorem~\ref{maingen}, $(\mathcal{F}_1, \dots, \mathcal{F}_k)$ has the strict cross-$t$-star property. Let $T = \{\{i\} \colon i \in [t]\}$. For each $i \in [k]$, we have $|\mathcal{F}_i(T)| = s_{n_i-t,r_i-t}$, and hence $\mathcal{F}_i(T) \in {\rm L}(\mathcal{F}_i,t)$ by Lemma~\ref{fsplemma2}. The result follows by Proposition~\ref{prop2}.~\hfill{$\Box$}

\end{document}